\documentclass[10pt,twocolumn,twoside]{IEEEtran}

\usepackage{stfloats}
\usepackage{caption}
\usepackage{bm}   
\usepackage{amsmath}
\usepackage{amssymb}
\usepackage{cuted}
\usepackage{indentfirst}  
\usepackage{graphicx}  
\usepackage{subfigure}
\usepackage{multirow} 
\usepackage{booktabs} 
\usepackage{cite}  
\usepackage{threeparttable}   
\usepackage{algorithm}
\usepackage{algorithmic}
\usepackage{arydshln}
\usepackage{color}
\usepackage[american]{babel}
\usepackage{microtype}
\usepackage[numbers,sort&compress]{natbib}  
\usepackage{steinmetz}
\usepackage{ifthen}
\usepackage{cases}
\usepackage{enumerate}
\usepackage{epstopdf}
\usepackage{float}
\usepackage{epstopdf}

\usepackage[cp1251]{inputenc}

\def\mL{\mathcal{L}}

\def\mF{\mathcal{F}}

\def\mbR{\mathbb{R}}

\def\mbS{\mathbb{S}}

\def\mL{\mathcal{L}}

\begin{document}
\title{An Analytical Formula for Stability Sensitivity Using SDP Dual}

\author{Jun~Wang,~\IEEEmembership{Graduate Student Member,~IEEE,}
        Yue~Song,~\IEEEmembership{Member,~IEEE,}
        David~J.~Hill,~\IEEEmembership{Life Fellow,~IEEE}
        and Yunhe~Hou,~\IEEEmembership{Senior Member,~IEEE}

\thanks{J. Wang, Y. Song, D. J. Hill and Y. Hou are with the Department of Electrical and Electronic Engineering, The University of Hong Kong, Hong Kong (e-mail: wangjun@eee.hku.hk; yuesong@eee.hku.hk; dhill@eee.hku.hk; yhhou@eee.hku.hk), D. J. Hill is also with the School of Electrical Engineering and Telecommunications, The University of NSW, Sydney (david.hill1@unsw.edu.au).}
\vspace{-20pt}
}

\markboth{IEEE TRANSACTIONS ON POWER SYSTEMS}  
{Wang \MakeLowercase{\textit{et al.}}: An Analytical Formula for Stability Sensitivity Using SDP Dual}

\maketitle

\begin{abstract}
In this letter, we analytically investigate the sensitivity of stability index to its dependent variables in general power systems. Firstly, we give a small-signal model, the stability index is defined as the solution to a semidefinite program (SDP) based on the related Lyapunov equation. In case of stability, the stability index also characterizes the convergence rate of the system after disturbances. Then, by leveraging the duality of SDP, we deduce an analytical formula of the stability sensitivity to any entries of the system Jacobian matrix in terms of the SDP primal and dual variables. Unlike the traditional numerical perturbation method, the proposed sensitivity evaluation method is more accurate with a much lower computational burden. This letter applies a modified microgrid for comparative case studies. The results reveal the significant improvements on the accuracy and computational efficiency of stability sensitivity evaluation.
\end{abstract}


%
\IEEEpeerreviewmaketitle
\vspace{-5pt}
\section{Introduction}
With an increasing penetrations of renewable energy sources appear in modern power systems, the stability problems with the dynamical behaviors significantly impact the system security \cite{b1}. Apart from controllers integrated in power systems, the re-dispatch may bring supplementary measures to enhance the system stability \cite{b13}. The stability sensitivities to the re-dispatch terms are of importance to determining the amount of re-dispatch, however, the stability sensitivity rarely has analytical expressions. The mainstream method to circumvent this obstacle is applying numerical perturbation. This method obtains the coefficient of the first-order Taylor expansion to formulate a linear approximation of stability constraints \cite{b13}-\cite{b42}. However, this numerical approach is CPU-consuming and inaccurate. By contrast, the analytical approach can find an accurate formula for stability sensitivity and is less computationally costly.

Practically, there are two main methods to analyze the small signal stability. One is the eigenvalue analysis, which uses the largest real part of eigenvalues of Jacobian matrix to assess the stability. The other applies Lyapunov equation, which is applied in this paper. \cite{b7} proposed the numerical perturbation-based sensitivity of the former. However, there exists no sensitivity formulae for the latter so far. Actually, unlike the eigenvalue analysis, the latter can not only check the stability, but characterize the convergence rate of post-disturbance oscillation. Consequently, a systematic study of analytical stability sensitivity based on Lyapunov equation needs to be established.

To fill the aforementioned research gap, this letter appropriately designs a analytical formula for stability sensitivity in a general power system. The stability index is described by a semi-definite program (SDP) based on Lyapunov equation. Because of the convexity and strong duality of the SDP, the analytical formula for stability sensitivity is deduced. The simulation results verify the accuracy and computational efficiency enhancements of the proposed sensitivity compared to numerical approaches. Hence, the proposed analytical stability sensitivity has widely applicability in planning and operating problems.

\section{Problem Formulation}

\subsection{Small Signal Stability Analysis of General Power Systems}
To address the stability issue, we first carry out a general small-signal model around an equilibrium point of a power system in the matrix form,
\begin{eqnarray} \label{SDM}
\left[ \begin{matrix}
    \Delta \dot{\bm{x}}\\\textbf{0}
\end{matrix} \right]
=\left[ \begin{matrix}
	\bm{A}&		\bm{B}\\
	\bm{C}&		\bm{D}\\
\end{matrix} \right]
\left[ \begin{matrix}
	\Delta  \bm{x}\\
	\Delta  \bm{y}\\
\end{matrix} \right]
\end{eqnarray}
where $\bm{x}$ denotes the vector of state variables and $\bm{y}$ denotes the vector of algebraic variables. The sub-matrices $\bm{A}$ and $\bm{B}$ relate to differential equations w.r.t. state variables and algebraic variables, respectively; and $\bm{C}$ and $\bm{D}$ relate to algebraic equations w.r.t. state variables and algebraic variables, respectively. It is common that matrix $\bm{D}$ is nonsingular. Then the system Jacobian matrix $\bm{J}$ can be obtained by eliminating $\Delta \bm{y}$
\vspace{-5pt}
\begin{equation}
{\Delta \dot{\bm{x}}} = (\bm{A} - \bm{B}{\bm{D}^{ - 1}}\bm{C})\bm{x}={\bm{J}}\Delta \bm{x}.
\end{equation}

Recalling the Lyapunov equation \cite{b34}, the system is asymptotically stable if and only if the existence of a symmetrical positive-definite real matrix $\bm{\Phi}$ such that
\vspace{-5pt}
\begin{equation}\label{stacondition1}
-{\bm{J}^T}\bm{\Phi} - \bm{\Phi}\bm{J} + \xi \bm{I}= \bm{0}
\end{equation}
with the relevant Lyapunov function $ \mL = \Delta {\bm{x}^T}\bm{\Phi}\Delta \bm{x}$, where $\xi$ is any given negative real number. To quantify the dynamic behavior, we further design the stability index $\eta$ as the solution to the following SDP problem
\begin{subequations}\label{stabilityindexdef}
\begin{align}
   \min &\quad {\eta}\\
   \textup{s.t.} &\quad - {\bm{J}}^T{\bm{\Phi}} - {\bm{\Phi}}{\bm{J}}+ {\eta}\bm{I} \succeq 0\label{orist}\\
   &\quad {\bm{\Phi}} - \epsilon \bm{I}\succeq 0\label{existence}\\
   &\quad  -{\bm{\Phi}} + \bm{I}\succeq 0\label{uni}
   \vspace{-10pt}
\end{align}
\end{subequations}
where the second constraint \eqref{existence} ensures the existence of the extreme points in the feasible region with $\epsilon$ being a very small positive number. Furthermore, we deduce \eqref{uni} to prevent the objective value of this problem from being infinity. If \eqref{uni} is not deduced, then assume $\bm{\Phi}$ is a feasible solution for this model, for any positive number $k$, $(k\bm{\Phi}, k\eta)$ must be another feasible solution. Consequently, the objective value may approach to infinity. If the system is stable, the stability index $\eta$ will be negative.

Based on the above definition of the stability index, the absolute value of stability index $|\eta|$ is the lower bound of post-disturbance convergence rate, since (4b) and (4d) give
\begin{equation}
\dot{\mL}=\bm{x}^T(\bm{J}^T\bm{\Phi}+\bm{\Phi}\bm{J})\bm{x}
\preceq {\eta}\bm{x}^T\bm{I}\bm{x} \preceq \eta\bm{x}^T\bm{\Phi}\bm{x}=\eta\mL.
\end{equation}

\vspace{-15pt}

\subsection{Analytical Stability Sensitivity Analysis}
The stability constraint is widely used for regulating the system stability,
\begin{equation}
\eta(\bm{d})<\bar{\eta}
\end{equation}
where $\bm{d}$ denotes a vector of controllable variables that can be used to enhance the stability index $\eta$, $\bar{\eta}$ represents the threshold. To enforce $\eta$ to satisfy (6), $\bm{d}$ should be properly adjusted by calculating the corresponding sensitivity. Due to the implicit function between $\eta$ and $\bm{d}$, the sensitivity is challenging to obtain. The mainstream sensitivity analysis is the numerical perturbation approach \cite{b14}, the stability sensitivity w.r.t. each element of $\bm{d}$ is estimated by:
\begin{equation}
\frac{\partial \eta}{\partial {d_i}}\bigg|_{{d_i}={d}_i^0}\approx \frac{\eta({d}_i^0+ \epsilon_p)-\eta({d}_i^0)}{\epsilon_p}
\end{equation}
which is inaccurate due to the evaluation result strongly depending on the value of perturbation $\epsilon_p$. Apart from the accuracy, the computational burden is high since the stability index needs to be calculated twice to just obtain the sensitivity w.r.t. a single variable $d_i$. Instead of adopting this numerical approach, this letter proposes an analytical formula to accurately calculate the stability sensitivity. From the viewpoint of the convexity of the SDP, the analytical sensitivity of $\eta$ to the system Jacobian entries $J_{ij}$ can be deduced by following steps.

For the convenience of our analysis, the original expression of \eqref{stabilityindexdef} needs to be rewritten as a standard parametric SDP,
\begin{equation}
\begin{split}
\min & \quad {\bm{c}^T}\bm{\phi}\\\
\textup{s.t.} & \quad\bm{F}(\bm{\phi},J_{ij}) \succeq 0\\
& \quad{\bm{F_i}} \in {\mbS^{3n}},i = 0,1,...,m
\end{split}
\label{primalSDP}
\end{equation}
where $\bm{c} = {(1,0,..,0)^T}$, vector $\bm{\phi}$ contains the stability index as well as entries of matrix $\bm{\Phi}$, $\bm{\phi}= {({\eta},{\Phi_{11}},...,{\Phi_{nn}})^T}$, $\bm{c},\bm{\phi} \in \mbR ^m, m = 1 + \frac{{n(n + 1)}}{2}$. By accounting constraints from \eqref{orist} to \eqref{uni}, $\bm{F}(\bm{\phi},J_{ij}) = {\bm{F_0}} + {\bm{F_1}}{\phi_1} + ... + {\bm{F_m}}{\phi_m}$, where each $\bm{F_i} = diag\{\bm{F_i}^{(1)},\bm{F_i}^{(2)},\bm{F_i}^{(3)}\}$, $\bm{F_i}^{(1)},\bm{F_i}^{(2)},\bm{F_i}^{(3)} \in {\mbS^{n}}$. $\mbS^{n}$ denotes the whole set of $n$-dimensional real symmetric matrices. More specifically, $\bm{F_0}^{(1)} = \bm{0},\bm{F_1}^{(1)} = \bm{I},\bm{F_i}^{(1)} =  - \bm{J}^T\bm{T_i} - \bm{T_i}\bm{J}$; $\bm{F_0}^{(2)} =  - \epsilon \bm{I},\bm{F_1}^{(2)} = \bm{0},\bm{F_i}^{(2)} = \bm{T_i}$; $\bm{F_0}^{(3)} = \bm{I},\bm{F_1}^{(3)} = \bm{0},\bm{F_i}^{(3)} =  - \bm{T_i}$, $\{\bm{T_i}\}$ denotes as a basis of $n$-dimensional symmetric matrices, $i = 2,...,m$.

Following the idea in \cite{b39}, the dual problem of the primal SDP problem \eqref{primalSDP} takes the following form
\begin{equation}
\begin{split}
\max & \quad- \textup{Tr}{\bm{F_0}}\bm{\Upsilon}\\
\textup{s.t.} & \quad \textup{Tr}{\bm{F}_i}\bm{\Upsilon} = {c_i},\\
& \quad\bm{\Upsilon} \succeq 0,\\
& \quad\bm{\Upsilon} \in {\mbS^{3n}},i = 1,...,m
\end{split}
\label{dualproblem}
\end{equation}
where $\bm{\Upsilon}$ is the corresponding dual variable. Due the convexity of the SDP problem, the duality gap at the optimum $(\bm{\phi},\bm{\Upsilon})$ between the original problem and the dual problem is zero
\begin{equation}
\textup{Tr}\bm{F}(\bm{\phi},J_{ij})\bm{\Upsilon}=0.
\end{equation}

Since $\bm{F}(\bm{\phi},J_{ij})\succeq 0$ and $\bm{\Upsilon}\succeq 0$, they together imply
\begin{equation}\label{dualitygap}
\frac{1}{2}\left[\bm{F}(\bm{\phi},J_{ij})\bm{\Upsilon} + \bm{\Upsilon}\bm{F}(\bm{\phi},J_{ij})\right] = \bm{0}
\end{equation}
\begin{equation}
\textup{Tr}{\bm{F}_i(J_{ij})}\bm{\Upsilon} = {c_i},i=1,2,...,m
\label{26}
\end{equation}
where \eqref{26} rewrites a constraint in \eqref{dualproblem}.

In order to simplify the presentation, we define a function $\textup{svec} : \mbS^n \rightarrow \mbR^{(m-1)}$ by $\bm{\upsilon} = \textup{svec}(\bm{\Upsilon}) \triangleq {[{\Upsilon_{11}},\sqrt 2 {\Upsilon_{12}},...,\sqrt 2 {\Upsilon_{1n}},{\Upsilon_{22}},\sqrt 2 {\Upsilon_{23}},...,{\Upsilon_{nn}}]^T}$. Then the $\textup{svec}$ function induces an isomorphism between these two vector spaces with a inner product
\begin{equation}
\textup{Tr}{\bm{F}_i(J_{ij})}\bm{\Upsilon}=\textup{svec}{({\bm{F}_i}({J_{ij}}))^T}\textup{svec}(\bm{\Upsilon})
\end{equation}
\begin{equation}\label{G1}
\textup{Tr}{\bm{F}_i(J_{ij})}\bm{\Upsilon}-c_i=\textup{svec}{({\bm{F}_i}({J_{ij}}))^T}\textup{svec}(\bm{\Upsilon})-c_i=0.
\end{equation}

Denote $\mF(J_{ij})\triangleq[\textup{svec}(\bm{F}_1({J_{ij}})),...,\bm{F}_m({J_{ij}}))]$, then \eqref{G1} can be simply formulated as
\begin{equation}\label{G2}
\mF(J_{ij})^T\bm{\upsilon}-c_i=0.
\end{equation}

Based on the duality gap being zero we will obtain
\begin{equation}\label{dualitysvec}
\textup{svec}(\frac{1}{2}\bm{F}(\bm{\phi},J_{ij})\bm{\Upsilon} + \bm{\Upsilon}\bm{F}(\bm{\phi},J_{ij})) = \bm{0}.
\end{equation}

For the simplicity of expression, the operator $\circledast$ is defined. For any three n-dimensional matrices $\bm{M}, \bm{N}, \bm{X}$, they satisfy the following equality
\begin{equation}
(\bm{M} \circledast \bm{N}) \triangleq \textup{svec}(\frac{1}{2}(\bm{N}\bm{X}\bm{M}^T+\bm{M}\bm{X}\bm{N}^T)),
\end{equation}
then the $\textup{svec}$ function corresponding to the duality gap \eqref{dualitysvec} can be rewritten as
\begin{equation}\label{G3}
(\bm{\Upsilon} \circledast \bm{I})\textup{svec}(\bm{F}(\bm{\phi},J_{ij}))= \bm{0}.
\end{equation}

By combining the results from the equalities \eqref{G2} and \eqref{G3}, a function $G$ which contains the original vector $\bm{\phi}$, dual variable $\textup{sevc}(\bm{\Upsilon})$ and parameter $J_{ij}$ is defined as
\begin{eqnarray} \label{G}
G(\bm{\psi},J_{ij})=\left[ \begin{matrix}
    G_1\\G_2
\end{matrix} \right]\triangleq
\left[ \begin{matrix}
    \mF(J_{ij})^T\bm{\upsilon}-c_i\\(\bm{\Upsilon} \circledast \bm{I})\textup{svec}(\bm{F}(x,J_{ij}))
\end{matrix} \right]
\end{eqnarray}
where $\bm{\psi}=[\bm{\phi}^T,\bm{\upsilon}^T]^T$ and $G(\bm{\psi},J_{ij})=\bm{0}$ at the optimum $\bm{\psi}^*$ for a given parameter $J_{ij}$.

Recalling the \textit{Implicit function theorem}, the parametric sensitivity $\frac{\partial \bm{\psi}}{\partial J_{ij}}$ in the original SDP \eqref{primalSDP} is obtained as
\begin{eqnarray}\label{parasensitivity}
{\left. {\frac{{\partial \bm{\psi}}}{{\partial {J_{ij}}}}} \right|_{{J_{ij}} = {J_{ij}^*}}} =  - G'{({\bm{\psi}^*},{J_{ij}^*})^{ - 1}} \cdot {\left. {\frac{{\partial G}}{{\partial {J_{ij}}}}} \right|_{{\bm{\psi}^*},{J_{ij}^*}}}
\end{eqnarray}
where the Jacobian matrix
\setlength{\arraycolsep}{1.2pt}
\begin{eqnarray}
G'({\bm{\psi}^*},{J_{ij}^*}) = {\left. {\left[ {\begin{array}{*{20}{c}}
0&{\mF {{({J_{ij}})}^T}}\\
{(\bm{\Upsilon} \circledast \bm{I})\mF ({J_{ij}})}&{\bm{F}(\bm{\phi},{J_{ij}}) \circledast \bm{I}}
\end{array}} \right]} \right|_{{\bm{\psi}^*},{J_{ij}^*}}}
\end{eqnarray}
and
\begin{eqnarray}
{\frac{{\partial G}}{{\partial {J_{ij}}}}} = {\left[ {\begin{array}{*{20}{c}}
{{{\left( {\frac{{\partial \mF}}{{\partial {J_{ij}}}}} \right)}^T}\bm{\upsilon}}\\
{(\bm{\Upsilon} \circledast \bm{I})(\frac{{\partial \mF}}{{\partial {J_{ij}}}}\bm{\phi} + \textup{svec}(\frac{{\partial {\bm{F}_0}}}{{\partial {J_{ij}}}}))}
\end{array}} \right]}.
\end{eqnarray}

By applying the chain rule, the sensitivity of stability index $\eta$ to the controllable variables $\bm{d}$ can be explicitly formulated,
\begin{equation}\label{anasensitivity}
\frac{\partial \eta}{\partial \bm{d}} = \sum_{i,j}\frac{\partial \eta}{\partial J_{ij}}\cdot \frac{\partial J_{ij}}{\partial \bm{d}}
\end{equation}
where the partial derivative $\frac{\partial J_{ij}}{\partial \bm{d}} $ can be easily obtained, because each Jacobian matrix entry is always explicit function w.r.t. $\bm{d}$.

Compared to the numerical perturbation method, the proposed analytical approach accurately provides the sensitivity with a much lower computational complexity. That is because the sensitivity is formulated by dual variables, which are by-products in the solution process and can be obtained without additional computations. Note that stability sensitivity information has vital importance and widespread applications in the planning and operation problems in a general power system. In the case study, we will apply a microgrid as an instance to demonstrate the effectiveness.

\vspace{-5pt}
\section{Case study}
\subsection{Simulation Setting}
 We apply a microgrid test model as described in \cite{b42}. The original line parameters are same as MATPOWER package \cite{b12}. To mimic real microgrids, we reduce the power injection at each bus to 10\% of its original value. The dispatchable DGs at bus 1, 5, 11, 16, 23, 29 are with the interfaced inverters, whose frequency and voltage droop gains $K_{gp}=0.75$, $K_{gq}=7.8$. The batteries are at bus 17, 21, 24 with adjustable droop gains $K_{bp}$, $K_{bq}$. With applying the unchanged microgrid structure and linear power flow model, the entry of Jacobian matrix $J_{ij}$ is has the explicit expression w.r.t $\bm{K_{b}}=(K_{bp}, K_{bq})$, which further determines the stability index.

In this section, we apply Monte-Carlo simulations to generate $10^3$ different scenarios which are used for calculating the value of $\frac{\partial \eta}{\partial \bm{K_b}}$ and cumulative computational time. Moreover, to illustrate the significant improvements on the accuracy and computational efficiency, the following three benchmarking numerical sensitivities as shown in (7) are selected, whose perturbation values $\epsilon_p = 10^{-1}, 10^{-2}$ and $10^{-3}$, respectively.
\vspace{-5pt}
\subsection{Accuracy and Computational Efficiency Improvement}
For quantifying the accuracy, we design two accurate degrees $\alpha_p$ and $\alpha_q$. Due to the space limit, only the formula of $\alpha_p$ is presented. The definition of $\alpha_q $ is similar to $\alpha_q$ with substituting $K_{bp}$ to $K_{bq}$ in \eqref{ap}.
\begin{equation}
\alpha_p = \left(1-\frac{\sum_{i=1}^{1000}{\left (\frac{\partial \eta}{\partial K_{bp}}\right)^{\textup{a}}_i - \left (\frac{\partial \eta}{\partial K_{bp}}\right)^{\textup{n}}_i / \left (\frac{\partial \eta}{\partial K_{bp}}\right)^{\textup{a}}_i}}{1000}\right) \times 100\%
\label{ap}
\end{equation}
where $\left(\frac{\partial \eta}{\partial K_{bp}}\right)^{\textup{a}}_i$ and $\left(\frac{\partial \eta}{\partial K_{bp}}\right)^{\textup{n}}_i$ are the analytical sensitivity and numerical sensitivity in the $i^{th}$
scenario, respectively. The superscripts ``n'' and ``a'' represent the same meanings in Table.~\ref{ARMSsensitivity}. It reveals the accuracy improvement by employing the analytical stability sensitivity compared to numerical approaches.
\begin{table}[htbp]
\vspace{-5pt}
\caption{Accuracy improvement through analytical sensitivity}
\vspace{-10pt}
\begin{center}
    \label{ARMSsensitivity}
    \renewcommand\tabcolsep{3.8pt} 
    \begin{threeparttable} 
        \begin{tabular}{ccccc} \toprule
             & Analytical Sensitivity & $\epsilon_p = 10^{-1}$ & $\epsilon_p = 10^{-2}$ & $\epsilon_p = 10^{-3}$\\
            \hline
            $\alpha_p$ & 100\%& 30.5\% & 36.6\%& 42.5\%\\
            $\alpha_q$ & 100\%& 28.6\% & 41.3\%& 43.7\%\\
            $T_\textup{cpu}$\tnote{1} & 128& 1982& 2135& 2029\\
            $r_\textup{t}$\tnote{2} & / & 93.5 \%& 94\%& 93.7\%\\
            \bottomrule
        \end{tabular}
        \begin{tablenotes}
            \item [1] $T_\textup{cpu}$ denotes the cumulative CPU time in seconds.
            \item [2] $r_\textup{t}$ denotes the time reduction which is calculated by $r_\textup{t}= \frac{T_\textup{cpu}^{\textup{n}}-T_\textup{cpu}^{\textup{a}}}{T_\textup{cpu}^{\textup{n}}}$.
        \end{tablenotes}
    \end{threeparttable}
\end{center}
\vspace{-15pt}
\end{table}

Apart from the accuracy, the computational efficiency is also dramatically improved by adopting the analytical sensitivity. The cumulative CPU time $T_\textup{cpu}$ for solving different sensitivities in 1000 scenarios is provided in Table.~\ref{ARMSsensitivity}. The proposed sensitivity can reduce $T_\textup{cpu}$ by more than $93\%$.
\vspace{-5pt}
\section{Conclusion}
We have proposed an analytical formula for stability sensitivity in this letter. We describe the stability index through the SDP problem by using the Lyapunov equation. With the dual property of SDP, we establish the analytical formula for stability sensitivity. The simulation results on a droop-controlled microgrid reveal that the proposed stability sensitivity is more accurate with a much lower computational complexity compared to numerical perturbation based-sensitivities.


\ifCLASSOPTIONcaptionsoff
  \newpage
\fi

{\footnotesize
\bibliographystyle{IEEEtran}

}

\end{document}